\documentclass[11pt,english]{article}

\usepackage{amsfonts,amsmath,amsthm,amscd,amssymb,latexsym,verbatim}

    \usepackage[T2A]{fontenc} 
    \usepackage[cp1251]{inputenc}
    \usepackage[ukrainian,russian]{babel}

\usepackage[dvips]{graphicx}

\theoremstyle{plain}
\newtheorem{theorem}{Теорема}
\newtheorem{problem}{Задача}
\newtheorem{lemma}{Лема}

\theoremstyle{definition}
\newtheorem{example}{Приклад}

\newcommand{\keywords}{\textbf{Key words. }\medskip}
\newcommand{\subjclass}{\textbf{MSC 2010. }\medskip}
\renewcommand{\abstract}{\textbf{Abstract. }\medskip}

\numberwithin{equation}{section}

\addto\captionsrussian{}

\begin{document}

\title{Деякі зауваження про системи куль, \\ які створюють тінь в точці}

\author{Тетяна Осіпчук}



\date{}

\maketitle

{\bf Анотація.} В даній роботі розглядаються задачі, пов'язані з відшуканням мінімального числа системи куль, які створюють тінь у фіксованій точці багатовимірного евклідового простору $\mathbb{R}^n$. Тут вираз "набір куль створює тінь  в точці"\/ означає, що  кожна пряма, яка проходить через задану точку, перетинає хоча б одну кулю з набору. В роботі встановлено нові властивості системи неперетинних куль з цетрами на сфері, що створюють тінь в довільній фіксованій точці внутрішності сфери у просторі $\mathbb{R}^3$. А також, побудовано систему з  $n+1$ неперетинних куль з рівними радіусами у просторі $\mathbb{R}^n$,  $ n\ge 3$, які створюють тінь у фіксованій точці простору.


{\bf Ключові слова:} {задача про тінь, система куль, сфера, еліпсоїд обертання, область, багатовимірний дійсний евклідовий простір.}

{\bf Аннотация.} В данной работе рассматриваются задачи, связанные с отысканием минимального количества шаров, которые создают тень в фиксированной точке многомерного евклидова пространства $\mathbb{R}^n$. Здесь выражение "набор шаров создает тень в точке"\/ означает, что любая прямая, проходящая через заданную точку, пересекает хотя бы один шар из набора.
 В работе установлены новые свойства системы непересекающихся шаров с центрами на сфере, которые создают тень в произвольной фиксированной точке внутренности сферы в пространстве $\mathbb{R}^3$. А также, построено систему из   $n+1$ непересекающегося шара с равными радиусами в пространстве $\mathbb{R}^n$,  $ n\ge 3$, которые создают тень в фиксированной точке пространства.

{\bf Ключевые слова:} {задача о тени, система шаров, сфера, эллипсоид вращения, область, многомерное действительное евклидово пространство.}


\begin{abstract} Problems, related to the determination of the minimal number of balls that generate a shadow at a fixed point in the multi-dimensional Euclidean space $ \mathbb{R}^n $, are considered in present work. Here, the statement "a system of balls generate shadow at a point"\/ means that any line passing through the point intersects at least one ball of the system.  New properties of balls that generate shadow at any fixed point inside of a sphere in space $\mathbb{R}^3$ are established and a system of $n+1$ disjoint balls with equal radii in $\mathbb{R}^n$,   $ n\ge 3$, that generate shadow at fixed point of the space is constructed in the work.
 \end{abstract}




\section{Вступ}

У 1982 роцi Г. Худайбергановим \cite{Hud1} була поставлена задача про тінь.

Нехай $x$ фіксована точка у багатовимірному дійсному евклідовому просторі $\mathbb{R}^n$.  Скажемо, що кулі  $\{B_i: i=\overline{1,m}\}$, $m<\infty$,  в $\mathbb{R}^n$, які  не містять $x$,  створюють в цій точці тінь,  якщо довільна пряма, що проходить через точку $x$, перетинає хоча б одну кулю з набору. Тоді задача про тінь може бути сформульована наступним чином:  {\it знайти мінімальне число відкритих (замкнених) та попарно неперетинних куль у просторі $\mathbb{R}^n$ з центрами на сфері $S^{n-1}$ та радіусами меншими радіуса сфери, які не містять центр сфери та створюють в ньому тінь}.

Задачу про тінь в такому формулюванні будемо називати класичною. Тут і надалі, під сферою $S^{n-1}$ будемо розуміти множину всіх точок простору $\mathbb{R}^n$, які знаходяться на однаковій відстані від деякої фіксованої точки простору \cite{Roz1_1}.

Класична задача про тінь була розв'язана Г. Худайбергановим при $n=2$: було показано, що для кола на площинi достатньо двох кругiв \cite{Hud1}. Там же було зроблено припущення про те, що i для випадку при $n>2$ мiнiмальне число таких куль рiвне $n$. Вона також цiкава з точки зору опуклого аналiзу тим, що є частковим випадком питання про належнiсть точки узагальнено опуклiй оболонцi сiм'ї компактних множин \cite{Zel2}.

В \cite{Zel2} Ю. Зелiнський довiв, що для $n=3$ трьох куль не достатньо, разом з тим, чотири кулi вже будуть створювати тiнь в центрi сфери. Там же встановлено, що для загального випадку у просторі $\mathbb{R}^n$, для довільного $n\ge 3$, мінімальною кількістю є  $n+1$ куля. Таким чином, класична задача про тінь повністю розв'язана.

Розглянемо наступні задачі, близькі до класичної задачі про тінь.

\begin{problem}\label{problem1}
Знайти мінімальне число відкритих (замкнених) та попарно неперетинних куль у просторі $\mathbb{R}^n$ з центрами на сфері $S^{n-1}$ та радіусами меншими радіуса сфери, які не містять фіксовану точку всередині сфери та створюють в цій точці тінь.
\end{problem}

\begin{problem}\label{problem2} \textup{(\cite{Zel8}\textup)}
Знайти мінімальне число відкритих (замкнених) та попарно неперетинних куль з рівними радіусами у просторі $\mathbb{R}^n$, які не містять фіксовану точку простору та створюють в цій точці тінь.
\end{problem}

В даній роботі доводяться дві теореми, які  частково розв'язують посталені задачі.

\begin{theorem} \label{theor3}
Нехай $S^2 (r)$ сфера з центром в нулі та радіусрм $r$ у просторі $\mathbb{R}^3$.  Позначимо через $n(x)$ найменше число відкритих куль, що не перетинаються, з центрами на сфері $S^2 (r)$ і таких, що не містять фіксовану точку $x\in\mathbb{R}^3$ та  створюють в цій точці тінь. Тоді $n(x)=3$,  для кожної точки $x\in\mathbb{R}^3$ такої, що $0\le |x|\le\frac{7}{9}r$.
\end{theorem}

\begin{theorem} \label{theor5}
Нехай $n(x)$ найменше число відкритих (замкнених) та попарно неперетинних куль з однаковими радіусами і таких, що не містять фіксовану точку $x\in\mathbb{R}^n$, $n\ge 3$, та  створюють в цій точці тінь.  Тоді  $n(x)\le n+1$.
\end{theorem}

В \cite{Zel3}, \cite{Zel4}, \cite{Osi6}  зроблено огляд цілої низки результатів, аналогічних  до класичної задачі про тінь, та їх узагальнень, отриманих Ю.Б.~Зелінським та його учнями.

У наступному розділі зроблено огляд тих результатів, які також частково дають розв'язок  задач \ref{problem1}, \ref{problem2}, та тих, які розв'язують деякі інші задачі про тінь.  Ці результати будуть сформульовані як леми, оскільки в межах даної роботи вони є допоміжними та використовуються для доведення теорем \ref{theor3}, \ref{theor5}.

\section{Допоміжні результати}

Наступний приклад дає один із способів побудови системи з  $n+1$   кулі  із задачі \ref{problem1}, які створюють тінь в центрі сфери.

\begin{example} \label{examp1} \textup{\cite{Zel2}}  Якщо у сферу вписати правильний $n$-вимірний симплекс
(див. \cite{Roz1_1}) та розмістити у його вершинах замкнені кулі з радіусами, величини яких дорівнюють половині довжини ребра симплекса, то ця система куль створить тінь в центрі сфери. Однак, кулі будуть попарно дотикатись одна одної, що протирічить умові задачі \ref{problem1}.  Нехай $a$ --- половина довжини ребра симплекса. Для досить малого $\varepsilon>0$ розглянемо систему куль, що складається з $n+1$ кулі, величини раіусів яких, відповідно, дорівнюють $a+\varepsilon$, $a-\varepsilon/2$, $a-\varepsilon/2^2$,  $a-\varepsilon/2^3$, $\ldots$, $a-\varepsilon/2^n$. Розмістимо ці кулі так, щоб вони дотикались одна одної, а їх центри утворювали симплекс, який, очевидно, незначно відрізняється від правильного. Тоді через центри цих куль проходить єдина сфера, в центрі якої відкриті кулі з тими ж радіусами створюють тінь. Якщо вихідні замкнені кулі трішки зменшити, то, в силу неперервності, такі кулі також будуть створювати тінь в центрі сфери.
\end{example}

 У \cite{Osi5,Osi6} розглядаються задачі про тінь для куль у просторах $\mathbb{R}^2$, $\mathbb{R}^3$  з центрами, розташованими на інших поверхнях.

\begin{lemma} \label{theor1} \textup{\cite{Osi5}}
Нехай задано видовжений еліпсоїд обертання з великою піввіссю $b$ та малою $a$ і нехай  $n(x_0)$ найменше число попарно неперетинних відкритих куль, з центрами на заданому еліпсоїді, які не містять його центр $x_0$ та створюють в $x_0$ тінь. Тоді,

1) $n(x_0)=3$, якщо $b/a>2\sqrt{2}$;

2) $n(x_0)>3$, якщо $b/a\le2\sqrt{2}$.
\end{lemma}

Наступна лема дає розв'язок задачі про тінь для кругів з центрами, розміщеними на довільній замкненій кривій на площині, а також  оцінку зверху мінімального числа куль, що створюють тінь в точці, з центрами на довільній замкненій поверхні у просторі  $\mathbb{R}^3$.

\begin{lemma} \label{theor2} \textup{\cite{Osi6}}
Нехай задано деяку обмежену область  $D\subset \mathbb{R}^3$  $\left( D\subset \mathbb{R}^2\right)$ і нехай $n(x)$ найменше число попарно неперетинних відкритих чи замкнених куль, з центрами на $\partial D$, які не містять фіксовану точку $x\in D$ та створюють в точці $x$ тінь. Тоді $n(x)\le 4$ $\left(n(x)=2\right)$.
\end{lemma}

Для її доведення застосовується наступна лема, яка в даній роботі буде використана також.

\begin{lemma} \label{corol1} \textup{\cite{Osi6}}
Нехай задано дві відкриті (замкнені) кулі $\{B_i=B(r_i)\}$, $i=1,2$, у просторі $\mathbb{R}^n$, які не перетинаються, з центрами на сфері $S^{n-1}(r)$ та радіусами $r_2\le r_1<r$. Тоді кожна куля, гомотетична кулі $B_2$  відносно центра сфери, з коефіцієнтом гомотетії $k_2$, не перетинає кожну кулю, гомотетичну кулі $B_1$  відносно центра сфери, з коефіцієнтом гомотетії $k_1$, якщо $k_2\ge k_1$.
\end{lemma}

В \cite{Zel7} розглядається задача про тінь для деякого набору куль з вільно розташованими центрами. Наступна лема дає оцінку знизу для числа неперетинних куль у просторі $\mathbb{R}^n$, які не містять фіксовану точку простору та створюють в ній тінь.

\begin{lemma} \label{theor4} \textup{\cite{Zel7}}
Нехай $n(x)$ найменше число відкритих (замкнених) та попарно неперетинних куль, які не містять фіксовану точку $x\in\mathbb{R}^n$, $n\ge 2$, та  створюють в цій точці тінь.  Тоді  $n(x)= n$.
\end{lemma}

В \cite{Zel8} ставиться задача про тінь для набору куль з вільно розташованими центрами, але з рівними радіусами і будується приклад з чотирьох відкритих (замкнених) та попарно неперетинних  куль $\{B_i\}$, $i=\overline{1,4}$,  однакового радіуса у просторі $\mathbb{R}^3$, які створюють тінь у фіксованій точці $x\in \mathbb{R}^3\setminus\cup_iB_i$. Таким чином встановлена наступна

\begin{lemma} \label{theor44} \textup{\cite{Zel8}}
Нехай $n(x)$ найменше число відкритих (замкнених) та попарно неперетинних куль з однаковими радіусами і таких, що не містять фіксовану точку $x\in\mathbb{R}^3$ та  створюють в цій точці тінь.  Тоді  $n(x)\le 4$.
\end{lemma}

Неважко показати, що для випадку простору $\mathbb{R}^2$, мінімальною  кількістю є дві кулі. В роботі \cite{Dak9} доведено, що ніякі три попарно неперетинні, замкнені (відкриті) кулі з рівними радіусами у просторі $\mathbb{R}^3$ не  створюють тінь у фіксованій точці простору по-за кулями. Таким чином, встановлено, що чотири є мінімальною кількістю вказаних куль у просторі $\mathbb{R}^3$.

\section{Доведення теорем 1 та  2}

У наступному прикладі запропоновано один із способів побудови  системи з трьох відкритих куль,  центри яких розташовано на  видовженому еліпсоїді обертання з відношенням великої півосі до малої $b/a>2\sqrt{2}$, і таких, що створюють тінь в центрі еліпсоїда.

\begin{example} \label{examp2}
Спочатку розглянемо еліпсоїд з відношенням $b/a=2\sqrt{2}$ та набір відкритих куль $B_i$, $i=\overline{1,3}$, заданих наступним чином.

Якщо центр першої кулі $B_1$ з радіусом, рівним малій півосі $a$, розмістити в основі цієї півосі, тоді відкритою залишиться тільки площина $\Sigma$, дотична до кулі $B_1$ в центрі еліпсоїда. Розглянемо кулі, дотичні до першої та з центрами на лінії обертання малої півосі. Неважко показати, що, якщо центр такої кулі прямує до точки, діаметрально протилежної до центу кулі $B_1$, тоді кут, який ця куля закриває для прямих, що проходять через центр еліпсоїда в площині  $\Sigma$, прямує до свого максимального значення $\pi/2$. Третя куля $B_3$, дотична до першої $B_1$ та з центром в основі великої півосі $b=2\sqrt{2} a$, також закриває кут
$$
\varphi=2\arcsin\left( \frac{\sqrt{a^2+(b)^2}-a}{b}\right)=\pi/2.
$$

\begin{center}
\includegraphics[width=10 cm]{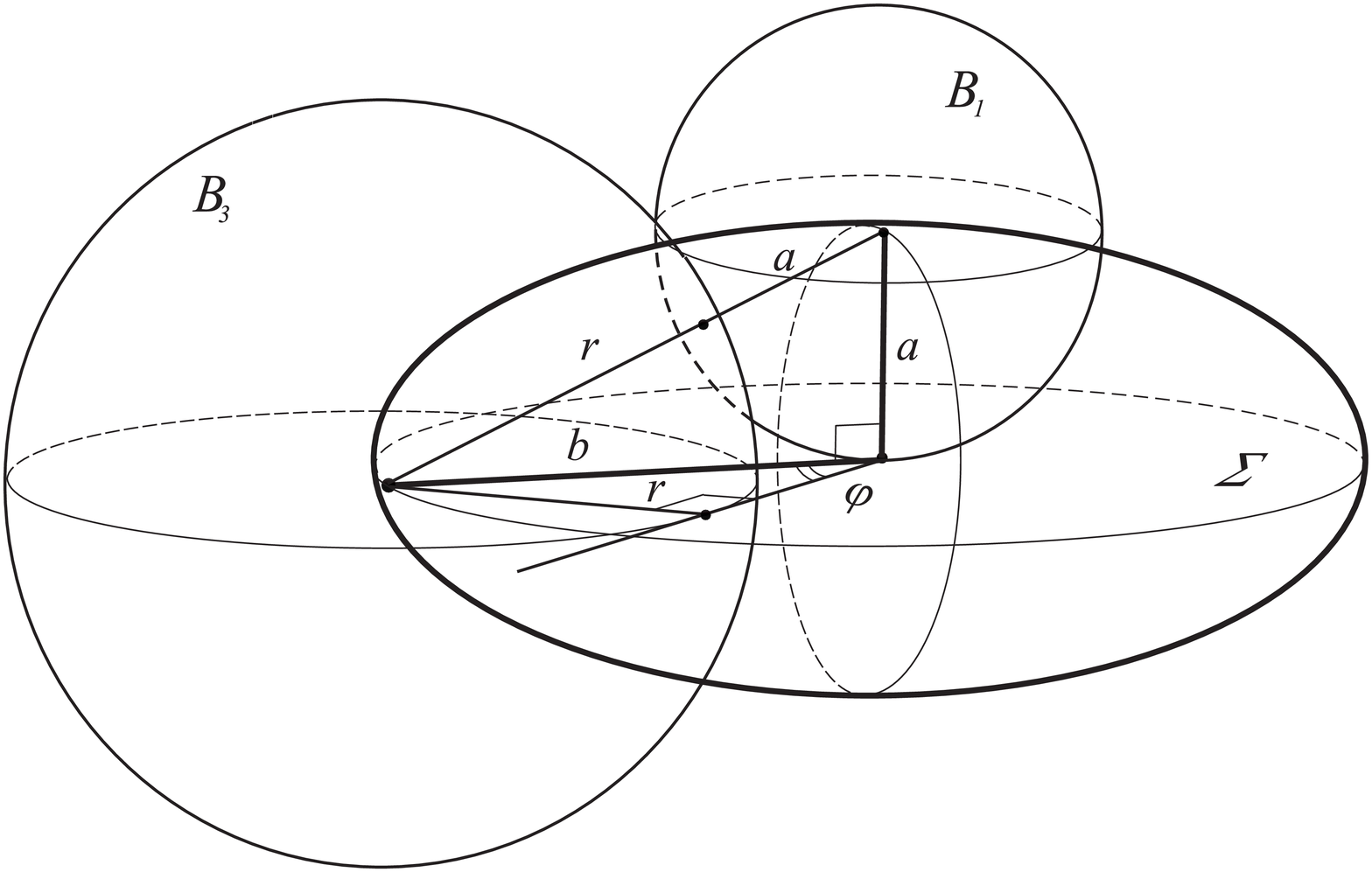}
\end{center}
\small\begin{center} Рис. 1 \end{center} \normalsize

Далі розглянемо еліпсоїд, центр та мала піввісь якого співпадають  з центром та малою піввіссю  попереднього еліпсоїда, а велика піввісь $b'>b$. Для таких еліпсоїдів  побудуємо систему відкритих куль $B_i'$, $i=\overline{1,3}$, наступним чином. $$B_1'\equiv B_1.$$

Куля $B_3'$ дотикається до кулі $B_1'$, а  її центр знаходиться в основі великої півосі $b'$. Тоді $B_3'$ закриває кут $\varphi(b')$ для прямих, що проходять через центр еліпсоїда в площині $\Sigma$, такий що
$$
\sin\frac{\varphi(b')}{2}=\frac{\sqrt{a^2+(b')^2}-a}{b'}.
$$
Оскільки
\begin{multline*}
\frac{d(\sin\varphi(b')/2)}{db'}=\frac{1}{\sqrt{a^2+(b')^2}}-\frac{\sqrt{a^2+(b')^2}-a}{(b')^2}= \\
=\frac{1}{\sqrt{a^2+(b')^2}}-\frac{\left(\sqrt{a^2+(b')^2}-a\right)\left(\sqrt{a^2+(b')^2}+a\right)}{(b')^2\left(\sqrt{a^2+(b')^2}+a\right)}=
\\=\frac{1}{\sqrt{a^2+(b')^2}}-\frac{1}{\sqrt{a^2+(b')^2}+a}>0
\end{multline*}
і $\varphi(b)=\pi/2$, тоді $\varphi(b')>\pi/2$ для $b'>b$.

Нарешті, центр кулі $B_2'$, яка дотикається до $B_1'$, розмістимо на лінії обертання малої півосі так, щоб кут, який вона закриває в площині  $\Sigma$ був більший  $\pi-\varphi(b')$.
\end{example}

\begin{proof}[{\bf Доведення теореми \ref{theor3}.}] Зафіксуємо точку $x$ внутрішності сфери $S^2 (x_0,r)$  на відстані $h>(7/9)r$ від її центра $x_0$, рис. 2. Побудуємо видовжений еліпсоїд обертання з центром в точці $x$, великою піввіссю $b=h+r$, розміщеною на прямій, що проходить через точку $x$ і центр сфери $x_0$, та малою піввіссю $a=\sqrt{r^2-h^2}$. Тоді неважко переконатись в тому, що для такого еліпсоїда відношення  $b/a>2\sqrt{2}$. Система із трьох куль, розміщених в точках перетину сфери з еліпсоїдом так, як це зроблено в прикладі \ref{examp2}, є шуканою, оскільки,  за лемою \ref{theor4}, не можливо побудувати систему куль, з кількістю менше трьох куль, які створюють тінь в довільній точці сфери.

\begin{center}
\includegraphics[width=8 cm]{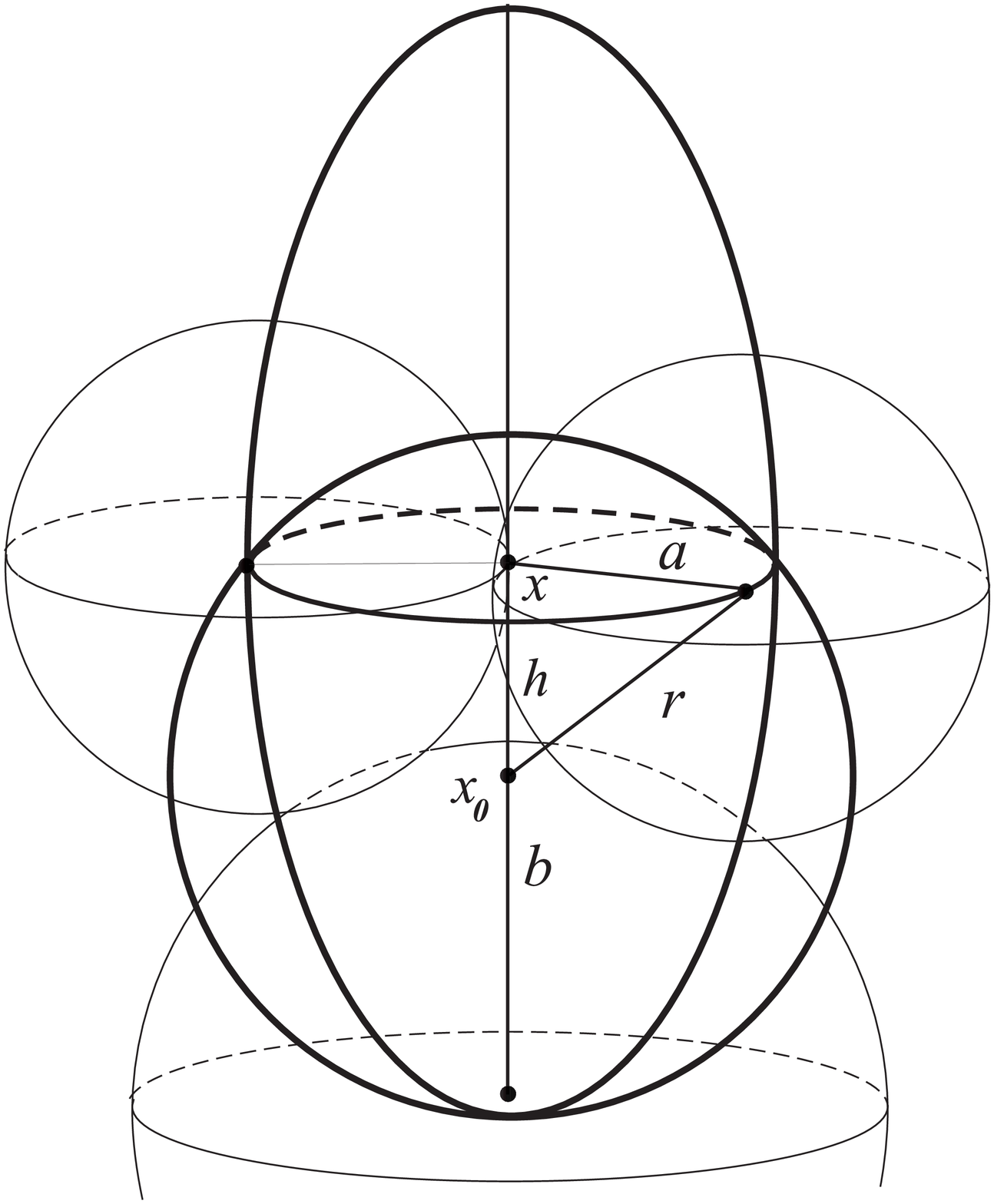}
\end{center}
\small\begin{center} Рис. 2 \end{center} \normalsize

\end{proof}

Для центра сфери задача про тінь розв'язана в \cite{Zel2}. Для решти точок всередині сфери $n(x)\le 4$ за лемою \ref{theor2}. Питання про те, чи для таких точок $n(x)=4$, залишається відкритим.

 \begin{proof}[{\bf Доведення теореми \ref{theor5}.}]
Зафіксуємо деяку точку $x\in\mathbb{R}^n$, $2<n<\infty$, та побудуємо систему з $n+1$ кулі $\{B_i=B(r_i)\}$, $i=1,\ldots,n+1$, з радіусами $r_i$, як у прикладі \ref{examp1}, розміщені на сфері з центром в точці $x$. Нехай, без обмеження загальності,  $r_1>\ldots>r_i>r_{i+1}>\ldots>r_{n+1}$.

До кожної кулі  $B_i$ застосуємо, відповідно, гомотетію з коефіцієнтом $k_i=r_1/r_i$. Тоді $k_{n+1}>\ldots>k_{i+1}>k_i>\ldots>k_1$ і отримана система складається з куль однакового радіуса, рівного $r_1$. За лемою \ref{corol1}, отримані кулі попарно не перетинаються і, за побудовою, не містять точку $x$ та створюють в ній тінь.
\end{proof}

Таким чином,  лему \ref{theor44} узагальнено на простір $\mathbb{R}^n$ довільної скінченної розмірності $n\ge 3$.

На даний момент залишається відкритим питання про те, чи справедливо, що в загальному випадку $n(x)=n+1$. Очевидно лише те, що згідно з лемою \ref{theor4}, це число не може бути меньше $n$.

\bigskip

CONTACT INFORMATION

\medskip
T.~Osipchuk\\01004 Ukraine, Kiev-4, 3, Tereschenkivska st., Institute of Mathematics of  NASU \\osipchuk.tania@gmail.com
\end{document}